\documentclass[12pt]{amsart}
\usepackage{amsfonts,amsmath,amsthm,amscd,amssymb,latexsym}
\usepackage{mathtools}
\usepackage{times}
\usepackage{amsmath,amsfonts,amssymb,amsopn,amscd,amsthm}
\usepackage{mathbbol}
\usepackage{mathrsfs}
\usepackage{graphicx,xspace}
\usepackage{epsfig}
\usepackage{enumerate} 
\usepackage{enumitem}
\usepackage{xcolor}
\usepackage[all]{xypic}

\usepackage[T1]{fontenc}
\usepackage[utf8]{inputenc}

\usepackage{hyperref}
\usepackage{psfrag}
\usepackage{bm}
\usepackage{youngtab}
\usepackage{tikz}
\usetikzlibrary{snakes}
\usepackage[all]{xy}
\usepackage{varioref}

\theoremstyle{plain}
\newtheorem{theoreme}{Theorem}[section]

\newtheorem{cor}[theoreme]{Corollary}

\theoremstyle{remark}

\newtheorem{ex}[theoreme]{Example}

\newtheorem{c-ex}[theoreme]{Counter-example}

\date{}

\DeclareUnicodeCharacter{03C7}{\chi}
\DeclareUnicodeCharacter{B2}{^{2}}

\DeclareMathOperator{\ty}{type}

\DeclareMathOperator{\Ind}{Ind}

\DeclareMathOperator{\diag}{diag}
\DeclareMathOperator{\Id}{Id}

\author[O. Tout]{Omar Tout}
\address{Department of Mathematics, Faculty of Sciences III, Lebanese University, Tripoli, Lebanon}
\email{omar.tout1988@gmail.com}

\title[Gelfand pairs involving finite abelian groups]{Gelfand pairs involving the wreath product of finite\\ abelian groups with symmetric groups}

 
\begin{document}
\maketitle

\begin{abstract} It is well known that the pair $(\mathcal{S}_n,\mathcal{S}_{n-1})$ is a Gelfand pair where $\mathcal{S}_n$ is the symmetric group on $n$ elements. In this paper, we prove that if $G$ is a finite group then $(G\wr \mathcal{S}_n, G\wr \mathcal{S}_{n-1}),$ where $G\wr \mathcal{S}_n$ is the wreath product of $G$ by $\mathcal{S}_n,$ is a Gelfand pair if and only if $G$ is abelian.
\end{abstract}

\section{Introduction}\label{sec1}

If $G$ is a finite group and $K$ is a subgroup of $G,$ we say that the pair $(G,K)$ is a \textit{Gelfand pair} if the algebra $C(G,K)$ of invariant complex valued functions on the double cosets of $K$ in $G$ is commutative. Equivalently, $(G,K)$ is a Gelfand pair if the algebra $\mathbb{C}[K\setminus G/K]$ spanned by the double cosets of $K$ in $G$ is commutative. If $(G,K)$ is a Gelfand pair then the algebra $C(G,K)$ has a particular basis whose elements are called zonal spherical functions of the pair $(G,K).$ 

Gelfand pairs are related to the theory of representation of finite groups. In fact, there is an equivalent definition for Gelfand pairs that uses induced representations. Explicitly, see \cite[chapter VII.1]{McDo}, $(G,K)$ is a Gelfand pair if and only if the induced representation of the trivial representation of $K$ to $G$ is multiplicity free. \\

Let $n$ be a positive integer and $\mathcal{S}_n$ be the symmetric group on $n$ elements. There are many Gelfand pairs involving symmetric groups in the literature. We present some of them along with references. First of all, the branching rule for the symmetric group shows that $(\mathcal{S}_n,\mathcal{S}_{n-1})$ is a Gelfand pair. For more details about this fact, the reader is invited to check \cite[Example 1.4.10]{ceccherini2010representation} and \cite{brender1976spherical}. The pair $(\mathcal{S}_n\times \mathcal{S}_n,\diag(\mathcal{S}_n)),$ where $\diag(\mathcal{S}_n)$ is the diagonal subgroup of $\mathcal{S}_n\times \mathcal{S}_n,$ is a Gelfand pair, see the Example $9$ in \cite[Section VII.1]{McDo}, and its zonal sphercial functions are the normalised irreducible characters of $\mathcal{S}_n.$ In \cite{strahov2007generalized}, Strahov shows that $(\mathcal{S}_n\times \mathcal{S}_{n-1},\diag(\mathcal{S}_{n-1}))$ is a Gelfand pair and its zonal spherical functions are the normalized generalized characters of $\mathcal{S}_n$ which form an orthogonal basis in the space of functions invariant with respect to conjugations by $\mathcal{S}_{n-1}.$ The pair $(\mathcal{S}_{2n}, \mathcal{H}_n),$ where $\mathcal{H}_n$ is the hyperoctahedral group $\mathcal{S}_2\wr \mathcal{S}_n,$ is a Gelfand pair studied in \cite[Section VII.2]{McDo} and \cite{toutejc}. Its zonal spherical functions appear in the development of zonal polynomials in terms of power-sums. The zonal polynomials are specialization of Jack polynomials which form a basis for the algebra of symmetric functions, see \cite{jack1970class} and \cite{jack1972xxv}. Recently, we showed in \cite{Tout2019OnTheSym} that $(\mathcal{H}_n\times \mathcal{H}_{n-1},\diag (\mathcal{H}_{n-1}))$ is a Gelfand pair. For more information about finite Gelfand pairs, the reader is invited to check \cite{scarabotti2008harmonic} and \cite{diaconis1988group}.\\

In \cite{aker2012parking}, Aker and Can prove that if $\Gamma$ is a finite abelian group then $(\widetilde{\Gamma^n} \ltimes \mathcal{S}_n, \mathcal{S}_n)$ is a Gelfand pair where $\Gamma^n=\Gamma\times \ldots \times \Gamma$ ($n$ copies), $\widetilde{\Gamma^n}$ is the quotient of $\Gamma^n$ by its diagonal subgroup $\lbrace (g,g,\ldots ,g)~~ |~~ g\in \Gamma\rbrace$ and $\widetilde{\Gamma^n} \ltimes \mathcal{S}_n$ is the semidirect product of $\widetilde{\Gamma^n}$ with $\mathcal{S}_n.$ In the particular case when $\Gamma=\mathbb{Z}_{r},$ the authors give explicit descriptions of the irreducible constituents of the multiplicity free representation $\Ind_{\mathcal{S}_n}^{\widetilde{\Gamma^n} \ltimes \mathcal{S}_n}(1).$\\

In \cite{stein2017littlewood}, Stein presents a complete list of irreducible representations of the wreath product $G\wr \mathcal{S}_n$ where $G$ is a finite group. In particular, he gives a branching rule for the group $G\wr \mathcal{S}_n$ similar to that of the symmetric group. Using this rule we will show the following main result of this paper 
\begin{theoreme}\label{th:main_result}
Let $G$ be a finite group then the pair $(G\wr \mathcal{S}_n,G\wr \mathcal{S}_{n-1})$ is a Gelfand pair if and only if $G$ is abelian.
\end{theoreme} 
We will only use the necessary definitions and facts from the representation theory of finite groups to prove this theorem. We eschew unuseful details and give references to the needed results instead of giving proofs. By doing so, we hope to present this result in a clear and short way.    

The paper is organised as follows. In Section \ref{sec2}, we review some useful tools from the theory of representation of finite groups and we define Gelfand pairs. In Section \ref{sec3}, we define multipartitions and we present the irreducible representations of the wreath product of a finite group $G$ by $\mathcal{S}_n.$ In Section \ref{sec4}, we prove Theorem \ref{th:main_result} the main result of this paper.

\section{Representation theory of finite groups and Gelfand pairs}\label{sec2}
In this paper we only consider representations over the complex numbers field $\mathbb{C}.$ We review in this section some definitions and well known results in the theory of representations of finite groups. We only present the necessery results that are useful to prove our main result. The reader is invited to check the book \cite{sagan2001symmetric} of Sagan for complete proofs. 

\subsection{Induced representations} Let $G$ be a finite group. A \textit{representation} of $G$ is a pair $(V,\rho)$ where $V$ is a finite dimensional vector space over $\mathbb{C}$ and $\rho :G \longrightarrow GL(V)$ is a group homomorphism. Here $GL(V)$ stends for the group of all invertible linear maps from $V$ to $V.$ The \textit{dimension} of a representation is the dimension of its vector space. For example, the trivial representation of $G$ is the one-dimensional representation $(\mathbb{C},1)$ where $1(g)(c)=c$ for any $g\in G$ and $c\in \mathbb{C}.$ 

Let $K$ be a subgroup of $G,$ $G/K=\lbrace g_1,\ldots ,g_r\rbrace$ be a set of representatives of the left cosets of $K$ in $G$ and $(V,\rho)$ be a representation of $K.$ Using $(V,\rho),$ we can make a representation $(\Ind_K^G V,\Ind_K^G \rho)$ of $G,$ called the \textit{induced representation} of $(V,\rho)$ to $G$ where 
$$\Ind_K^G V=\bigoplus_{i=1}^rg_iV$$
 with $g_iV$ being an isomorphic copy of $V$ whose elements are written as $g_iv$ with $v\in V.$ For each $g\in G$ and each $g_i$ there is an $h_i\in K$ and $j(i)\in \lbrace 1,\ldots, r\rbrace$ such that $g g_i = g_{j(i)} h_i$ and $\Ind_K^G \rho(g)$ is defined by:
$$\Ind_K^G \rho(g)\Big(\sum_{i=1}^rg_iv_i\Big)=\sum_{i=1}^rg_{j(i)} \rho(h_i)(v_i).$$
For example, in the particular case of the trivial representation of $K,$ the induced representation can be identified with the representation $(\mathbb{C}[G/K],1_K^G)$ of $G,$ where $\mathbb{C}[G/K]$ is the $r$-dimensional vector space spanned by the left cosets $g_iK$ and 
$$1_K^G(g)(k_iK)=(gk_i)K \text{ for any $g\in G$ and any representative $g_i$ of $G/K$}.$$  

\subsection{Facts from representation theory of finite groups} A vector subspace $W$ of $V$ is said to be a \textit{subrepresentation} of $(V,\rho)$ if $\rho(g)(w)\in W$ for any $w\in W.$ The representation $(V,\rho)$ is said to be \textit{reducible} if it has a subrepresentation other than $\lbrace 0\rbrace$ and $V.$ Otherwise it is said to be \textit{irreducible}. If $(V_1,\rho_1)$ and $(V_2,\rho_2)$ are two representations of $G$ then $(V_1\oplus V_2,\rho_1\oplus \rho_2),$ where $V_1\oplus V_2$ is the direct sum of $V_1$ and $V_2,$ is a representation of $G$ with:
$$(\rho_1\oplus \rho_2)(g)(v_1+v_2):=\rho_1(g)(v_1)+\rho_2(g)(v_2),$$
for any $g\in G, v_1\in V_1$ and $v_2\in V_2.$

A representation $(V,\rho)$ is equivalent to an action of $G$ on $V,$ that is $gv=\rho(g)(v)$ for any $g\in G$ and $v\in V.$ To simplify notations, it will be convenient for us to omit the homomorphism and say that $V$ is a representation of $G.$
If $V$ and $W$ are two representations of $G$ then $\theta:V\rightarrow W$ is said to be a \textit{representations homomorphism} (or an \textit{intertwiner}) if $\theta$ is an homomorphism of vector spaces that satisfies $\theta(gv)=g\theta(v)$ for any $g\in G$ and any $v\in V.$ Two representations are \textit{equivalent} if there exists a representations isomorphism between them. Otherwise, they are called \textit{inequivalent representations}. We will present below some well known facts in the theory of representation of finite groups. For complete proofs, the reader may refer to the book \cite{sagan2001symmetric} of Sagan.\\ 

The number of inequivalent irreducible representations of $G$ is equal to the number of conjugacy classes of $G,$ see \cite[Proposition $1.10.1$]{sagan2001symmetric}. Let $(V, \rho)$ and $(W, \sigma)$ be two irreducible representations. Then the Schur lemma (see \cite[Theorem $1.6.5$]{sagan2001symmetric}) states that any representation homomorphism $\theta : V \longrightarrow W$ is either the zero map (if $\rho\not\sim \sigma$) or a vector space isomorphism (if $\rho\sim \sigma$). As a consequence of the Schur lemma, if $(V, \rho)$ and $(W, \sigma)$ are two irreducible representations over the same
space, then any representation homomorphism $\theta : V \longrightarrow V$ is either the zero map (if $\rho\not\sim \sigma$) or $\theta=c\Id_V$ for a certain nonzero $c\in \mathbb{C},$ where $\Id_V$ is the identity map on $V$ (if $\rho\sim \sigma$).
Let $\lbrace V^1, \ldots, V^l\rbrace$ be a complete list of irreducible pairwise inequivalent representations of $G$ then by \cite[Proposition $1.10.1$]{sagan2001symmetric}
$$|G|=\sum_{i=1}^l (\dim V^i)^2.$$ 
In addition, by  Maschke's theorem (see \cite[Theorem $1.5.3$]{sagan2001symmetric}) any finite dimensional representation $V$ of $G$ can be decomposed as the direct sum of irreducible representations of $G$  
$$V\simeq \bigoplus_{i=1}^{l}m_iV^i,$$
where $m_i\in \mathbb{N}$ and $m_iV^i:=V^i\oplus V^i\oplus \cdots\oplus V^i$ ($m_i$ copies) for any $1\leq i\leq l.$ We say that the representation $V$ of $G$ is \textit{multiplicity free} if $m_i\leq 1$ for any $1\leq i\leq l.$

\subsection{Gelfand pairs} Let $(G,K)$ be a pair where $G$ is a finite group and $K$ is a subgroup of $G.$ Consider $C(G,K),$ the algebra of invariant complex valued functions on the double cosets of $K$ in $G,$ that is
$$C(G,K):=\lbrace f:G\longrightarrow \mathbb{C}~~ |~~ f(kxk')=f(x) \text{ for any $x\in G$ and any $k,k'\in K$} \rbrace.$$
The multiplication in $C(G,K)$ is defined by the convolution product of functions :
$$(fg)(x)=\sum_{y\in G}f(y)g(y^{-1}x)\text{ for any $f,g\in C(G,K)$}.$$

The pair $(G,K)$ is said to be a \textit{Gelfand pair} if its associated algebra $C(G,K)$ is commutative. Equivalently, $(G,K)$ is a Gelfand pair if the algebra $\mathbb{C}[K\setminus G/K]$ spanned by the double cosets of $K$ in $G$ is commutative. In \cite[$(1.1)$ page $389$]{McDo}, it is shown that the pair $(G,K)$ is a Gelfand pair if and only if the induced representation $\mathbb{C}[G/K]$ of the trivial representation of $K$ to $G$ is multiplicity free.

\section{Multipartitions and the irreducible representations of $G\wr \mathcal{S}_n$}\label{sec3}
In this section we present a complete list of inequivalent irreducible representations of $G\wr \mathcal{S}_n$ for a finite group $G.$ As in the previous section, we don't give the proofs of the presented results. The reader is invited to read the paper \cite{stein2017littlewood} of Stein for more details. 
 
\subsection{Partitions} A \textit{partition} $\lambda$ is a weakly decreasing list of positive integers $(\lambda^1,\ldots,\lambda^l)$ where $\lambda^1\geq \lambda^2\geq\ldots \geq\lambda^l\geq 1.$ The $\lambda^i$ are called the \textit{parts} of $\lambda$; the \textit{size} of $\lambda$, denoted by $|\lambda|$, is the sum of all of its parts. If $|\lambda|=n$, we say that $\lambda$ is a partition of $n.$ We will also use the exponential notation $\lambda=(1^{m_1(\lambda)},2^{m_2(\lambda)},3^{m_3(\lambda)},\ldots),$ where $m_i(\lambda)$ is the number of parts equal to $i$ in the partition $\lambda.$ In case there is no confusion, we will omit $\lambda$ from $m_i(\lambda)$ to simplify our notation. If $\lambda=(1^{m_1(\lambda)},2^{m_2(\lambda)},3^{m_3(\lambda)},\ldots,n^{m_n(\lambda)})$ is a partition of $n$ then $\sum_{i=1}^n im_i(\lambda)=n.$ We will dismiss $i^{m_i(\lambda)}$ from $\lambda$ when $m_i(\lambda)=0.$

If $\lambda=(\lambda^1,\ldots,\lambda^l)$ is a partition of $n,$ there are several natural ways to extend it to a partition of $n+1.$ This can be done by adding $1$ to $\lambda^1,$ adding a part $\lambda^{l+1}$ equal to one, adding $1$ to the parts $\lambda^i$ that satisfy $\lambda^{i-1}>\lambda^i\geq \lambda^{i+1}$ for $2\leq i\leq l-1$ or by adding $1$ to the parts $\lambda^{i+1}$ that satisfy $\lambda^{i-1}=\lambda^i> \lambda^{i+1}$ for $2\leq i\leq l-1.$ We will denote by $\lambda^+$ the set of all extended partitions from $\lambda.$ For example, $$(3,3,2,2,2,1)^+=\lbrace (4,3,2,2,2,1),(3,3,2,2,2,1,1),(3,3,3,2,2,1),(3,3,2,2,2,2)\rbrace.$$
\subsection{Branching rule for $\mathcal{S}_n$} 
It is well known that partitions of $n$ are used to index all of the conjugacy classes and the irreducible representations of $\mathcal{S}_n.$ The irreducible representation of $\mathcal{S}_n$ that corresponds to the partition $\lambda$ of $n$ is called \textit{Specht module} and usually denoted by $S^\lambda,$ see \cite[Theorem $2.4.6$]{sagan2001symmetric}. The branching rule shows that the inducing representation of an irreducible representation $\mathcal{S}^\lambda$ of $\mathcal{S}_{n-1}$ to $\mathcal{S}_n$ has a nice decomposition into irreducible representations of $\mathcal{S}_n.$ Explicitly, if $\lambda$ is a partition of $n-1$ then by \cite[Theorem $2.8.3$]{sagan2001symmetric} we have the following formula:

\begin{equation*}
\Ind_{\mathcal{S}_{n-1}}^{\mathcal{S}_n}S^\lambda=\bigoplus_{\delta\in \lambda^+}S^\delta.
\end{equation*} 

In particular, since the trivial representation of $\mathcal{S}_{n-1}$ corresponds to the Specht module $S^{(n-1)}$ by \cite[Example $2.3.6$]{sagan2001symmetric}, we have

\begin{equation*}
1_{\mathcal{S}_{n-1}}^{\mathcal{S}_n}=\Ind_{\mathcal{S}_{n-1}}^{\mathcal{S}_n}S^{(n-1)}=S^{(n)}\oplus S^{(n-1,1)}
\end{equation*} 
and $(\mathcal{S}_n,\mathcal{S}_{n-1})$ is a Gelfand pair. 
\subsection{Irreducible representations and a branchig rule for $G\wr \mathcal{S}_n$} Let $G$ be a finite group. The \textit{wreath product} $G \wr \mathcal{S}_n$ is the group with underlying set $G^n\times \mathcal{S}_n$ and product defined as follows:
$$((\sigma_1,\ldots ,\sigma_n); p).((\epsilon_1,\ldots ,\epsilon_n); q)=((\sigma_1\epsilon_{p^{-1}(1)},\ldots ,\sigma_n\epsilon_{p^{-1}(n)});pq),$$
for any $((\sigma_1,\ldots ,\sigma_n); p),((\epsilon_1,\ldots ,\epsilon_n); q)\in G^n\times \mathcal{S}_n.$ 

Suppose now that $\lbrace V^1, \ldots, V^l\rbrace$ is a complete list of irreducible pairwise inequivalent representations of $G.$ We assume that $V^1$ is the trivial representation of $G.$ An \textit{$l$-multipartition} of $n$ is a tuple $\Lambda=(\lambda_1,\lambda_2,\ldots ,\lambda_l)$ where $\lambda_i$ is a partition for each $1\leq i\leq l$ and 
\begin{equation*}
|\Lambda|:=\sum_{i=1}^l|\lambda_i|=n.
\end{equation*}
The irreducible representations of $G\wr \mathcal{S}_n$ are indexed by $l$-multipartitions of $n.$ We will denote by $S^{\Lambda}$ the irreducible representation of $G\wr \mathcal{S}_n$ that corresponds to an $l$-multipartition $\Lambda$ of $n.$ The set $\lbrace S^{\Lambda} ~~ |~~ \text{ $\Lambda$ is an $l$-multipartition of $n$} \rbrace$ forms a complete list of irreducible representations of $G\wr \mathcal{S}_n,$ see for example \cite{stein2017littlewood}. In \cite[Theorem $5.1$]{stein2017littlewood}, there is a branching rule for the group $G\wr \mathcal{S}_n$ similar to that of the symmetric group $\mathcal{S}_n$ presented in the previous subsection.  Explicitly, if $\Lambda=(\lambda_1,\lambda_2,\ldots ,\lambda_l)$ is an $l$-multipartition of $n-1$ then we have the following formula:

\begin{equation*}
\Ind_{G\wr \mathcal{S}_{n-1}}^{G\wr \mathcal{S}_n}S^\Lambda=\bigoplus_{i=1}^{l}\dim V^i\bigoplus_{\delta\in \lambda_i^+}S^{(\lambda_1,\ldots ,\lambda_{i-1},\delta,\lambda_{i+1},\ldots ,\lambda_l)}.
\end{equation*} 

In particular,

\begin{eqnarray}\label{eq:induced_representation}
1_{G\wr \mathcal{S}_{n-1}}^{G\wr \mathcal{S}_n}&=&\Ind_{G\wr \mathcal{S}_{n-1}}^{G\wr \mathcal{S}_n}S^{\big((n-1),\emptyset,\emptyset,\ldots ,\emptyset\big)}  \\
&=&S^{\big((n),\emptyset,\emptyset,\ldots ,\emptyset\big)}\oplus S^{\big((n-1,1),\emptyset,\emptyset,\ldots ,\emptyset\big)} \bigoplus_{i=2}^{l}\dim V^i. S^{\big((n-1),\emptyset,\ldots ,\emptyset,(1),\emptyset,\ldots \emptyset,\emptyset\big)}, \nonumber
\end{eqnarray} 
where $(1)$ is in the $i^{th}$ position of $S^{\big((n-1),\emptyset,\ldots ,\emptyset,(1),\emptyset,\ldots \emptyset,\emptyset\big)}.$

\section{proof of the main result}\label{sec4}

By Equation (\ref{eq:induced_representation}), the pair $(G\wr \mathcal{S}_n,G\wr \mathcal{S}_{n-1})$ is a Gelfand pair if and only if all the irreducible representations $V^i$ of $G$ are one dimensional. In the following paragraph we will show that all the irreducible representations of a finite group $G$ over the field of complex numbers $\mathbb{C}$ are one dimensional if and only if $G$ is abelian. Thus, we prove Theorem \ref{th:main_result}.

To start with, suppose that $G$ is a finite abelian group and that $(V,\rho)$ is an irreducible representation of $G.$ For any $g\in G,$ $\rho(g):V\rightarrow V$ is a $G$-homomorphism, that is $\rho(g)\circ \rho(h)=\rho(h)\circ \rho(g)$ for any $h\in G.$ This can be checked easily since $\rho(g)\circ \rho(h)=\rho(gh)=\rho(hg)=\rho(h)\circ \rho(g).$ But $(V,\rho)$ is irreducible, it follows then from Schur's lemma that for any $g\in G,$ $\rho(g)=c_g\Id_V,$ where $c_g\in \mathbb{C}.$ Let $v$ be a nonzero vector of $V,$ the subspace $\prec v\succ$ generated by $v$ is $G$-invariant since for any $av\in \prec v\succ$ with $a\in \mathbb{C}$ we have $\rho(g)(av)=a\rho(g)(v)=ac_gv\in \prec v\succ.$ Thus $\prec v\succ$ is a subrepresentation of $V$ which is irreducible. Therefore $\prec v\succ=V$ and $\dim_{\mathbb{C}}V=1.$ In the opposite direction, we have 
\begin{equation*}
|G|=\sum_{i=1}^l (\dim V^i)^2.
\end{equation*} 
Thus if all the $V^i$ are one dimensional then the number of irreducible representations of $G,$ which is $l,$ is equal to $|G|.$ Thus the number of conjugacy classes of $G$ is $|G|,$ that is the conjugacy class of any $g\in G$ is $\lbrace g\rbrace.$ This means that for any $g, h\in G,$ $hgh^{-1}=g$ or equivalently $hg=gh$ and $G$ is thus abelian.

\bibliographystyle{abbrv}
\bibliography{../biblio}

\end{document}